\documentclass[11pt,hyp,]{amsart}
\usepackage{hyperref}
\hypersetup{nesting=true,debug=true,naturalnames=true}
\usepackage{graphicx,amssymb,upref}

\newtheorem{theorem}{Theorem}
\newtheorem{lemma}[theorem]{Lemma}

\title{Extended gambler's ruin problem}
\author{Theo van Uem}  
\address{Amsterdam University of Applied Sciences, Amsterdam, The Netherlands.} 
\email{tjvanuem@gmail.com}  

\begin{document}
\hbadness=99999


\begin{abstract}
In the extended gambler's ruin problem we can move one step forward or backward (classical gamblers ruin problem), we can stay where we are for a time unit (delayed action) or there can be absorption in the current state (game is terminated without reaching an absorbing barrier). We obtain absorption probabilities,  probabilities for maximum and minimum values of  the ruin problem, expected time until absorption and the value of the game. We also investigate asymptotic behavior of absorption probabilities and  expected time until absorption  We  introduce  a conjugate version of our random walk.
\end{abstract}
\maketitle

\section{Introduction}
The gambler's ruin problem is a special random walk. Random walk can be used in various disciplines: in physics as a simplified model of Brownian motion, in ecology to describe individual animal movements and population dynamics, in statistics to analyze sequential test procedures, in economics to model share prices and their derivatives, in medicine and biology where absorbing barriers give a natural model for a wide variety of phenomena.
In Feller \cite{fel} there is a complete chapter (XIV) devoted to random walk and ruin problems. El-Shehawey et al. \cite{els} consider  gambler's ruin problem in the case that the probabilities of winning/losing a particular game depend on the amount of the current fortune with ties allowed. Yamamoto \cite{yam}, using hypergeometric functions,  treats a random walk which  moves either rightwards or leftwards, and in addition introduces the ‘halt’: the walker does not move for a time unit.
In this paper we investigate an extended one dimensional random walk. We call it a
 $[pqrs]$ walk, where $p$ is the one-step forward probability, $q$ one-step backward, $r$ the probability to stay for a time unit in the same position and $s$ is the probability of absorption in the current state $(p+q+r+s=1,\ pqs>0)$. We also use absorbing barriers to model the ruin problem. Proceeding along these lines, we shall model and analyze several extended ruin problems:\\
 a. The  gambler's ruin problem on $[0,N]$ with absorbing barriers in $0$ and $N$  and  actions in each non-barrier state: moving one step forward or backward, stay for a moment in the same position or  terminate the process (absorption in current state has occurred).\\
 b. A ruin problem on $[0,\infty)$, where $0$ is the only  absorbing barrier and in all other states we can move one step forward or backward, we can stay for a moment or absorption in the current state occurs.\\ 
 c. Unlimited resources: a ruin problem on $(-\infty,\infty)$, starting in state $0$ (gain= $0$) and moving to left or right or staying for a moment or absorption occurs and the game is terminated. There are no absorbing barriers in this case.\\
  In section \ref{s2} we solve a set of difference equations which is fundamental in our paper.
Section \ref{s3} deals with absorption probabilities, including asymptotic behavior.
Section \ref{s4} covers the expected time until absorption  in a $[pqrs]$ random walk, also including asymptotic behavior. 
In section \ref{s5} we obtain results for maximum and minimum of the random walk.
 Section \ref{s6} studies the value of the game. Section \ref{s7}  introduces a conjugate  random walk. 

\section{A related set of difference equations} \label{s2}

For a discrete Markov chain  we define the expected number of visits to state j when starting in state i by:
\begin{equation*}
x_j=x_{i,j}=\sum_{k=0}^{\infty}p_{i,j}^{(k)}
\end{equation*}
We start in state $i_0.$\\
Last step analysis gives:
$ x_ n  =\delta   (n,i_0 )+  px_{n- 1}+   qx_{ n+ 1}+ rx_n$.\\
The next theorem is fundamental for the rest of our paper.
\begin{theorem}
 \label{t1}
The set of difference equations:
\begin{equation} \label{1}
(1 - r) x_ n  =\delta   (n,i_0 )+  px_{n- 1}+   qx_{ n+ 1} \quad	(a<n<b )			
\end{equation}
where $ pq > 0, \ 	p + q+  r < 1,\  	$
has  solutions:			
\begin{equation} \label{2}
x_n=\left\{\begin{array}{l}\zeta \xi_1^{n-i_0}+C_1 \xi_1^{n}+C_2 \xi_2^{n}\ \ \ \ \ \ \ \    (a\leq n \leq i_0) \\
\zeta \xi_2^{n-i_0}+C_1 \xi_1^{n}+C_2 \xi_2^{n}\ \ \ \ \ \ \ \    ( i_0\leq n \leq b)\end{array}\right.
\end{equation}	
where:
\[
\xi_1=\frac{(1-r)+\sqrt{(1-r)^2-4pq}}{2q}>1
\]
\[
0<\xi_2=\frac{(1-r)-\sqrt{(1-r)^2-4pq}}{2q}<1
\]

			\[
			\zeta=	[(1 - r)^2-	4 pq	]^{-\frac{1}{2}}		
			\]	
\end{theorem}	
\begin{proof}
																				
General solution of homogeneous part of \eqref{1} is:						\[
x_n=C_1 \xi_1^{n}+C_2 \xi_2^{n} \ \ (n\in \mathbb Z)\]
 where $\xi_1$ and $\xi_2$ are the solutions of:
 \[q \xi^2-  (1 - r )\xi +   p =0
 	\]
		
A particular solution of \eqref{1} is (verified by substitution):						\begin{equation*} 
x_n=\frac{1}{2\pi}	\int_{-\pi}^{\pi}\frac{\exp [-\mathrm i\theta(n-i_0)]\mathrm{d}\theta}{(1-r)-p\exp(\mathrm i\theta)	-qexp(-\mathrm i\theta)}										
	\end{equation*}		
Substituting $z =e^{-\mathrm i\theta}$	gives
\[x_n=\frac{\mathrm i}{2\pi}\oint\frac{z^{n-i_0}\mathrm d z }{qz^2-(1-r)z+p}=
\frac{\mathrm i}{2\pi}\oint\frac{z^{n-i_0}\mathrm d z }{q(z-\xi_1)(z-\xi_2)}\]
where the integration is counterclockwise around the circle $|z|=1$.\\
After applying the residue theorem we obtain a particular solution:	
\begin{equation*} 
x_n=\left\{\begin{array}{l}\zeta \xi_1^{n-i_0} \hspace{3em}    (n \leq i_0) \\
\zeta \xi_2^{n-i_0}\hspace{3em}    (n \geq i_0)\end{array}\right.
\end{equation*}
General solution:
\begin{equation} \label{3}
x_n=\left\{\begin{array}{l}\zeta \xi_1^{n-i_0}+C_{11} \xi_1^{n-i_0}+C_{12} \xi_2^{n-i_0} \hspace{3em}    (n \leq i_0) \\
\zeta \xi_2^{n-i_0}+C_{21} \xi_1^{n-i_0}+C_{22} \xi_2^{n-i_0}\hspace{3em}    (n \geq i_0)\end{array}\right.
\end{equation}	
By substituting $n=i_0$ twice in \eqref{3} and taking $n=i_0$ in \eqref{1} we get:

\begin{equation*} 
\left\{\begin{array}{l}C_{11} +C_{12} = C_{21} +C_{22}\\
(1-r-q\xi_2)C_{11} +(1-r-q\xi_1)C_{12} = q\xi_1C_{21} +q\xi_2C_{22}
\end{array}\right.
\end{equation*}
We apply Cramer's rule, with $C_{11}$ and $C_{12}$ as variables, to obtain:

 $C_{11}=C_{21}$ and $C_{12}=C_{22}.$
 \end{proof}
The $x_n$  are unique: given an arbitrary solution of   \eqref{1}, the constants $C_1$ and $C_2$ can be chosen so that \eqref{2} will agree with it for two consecutive values of $n$. From these two $x_n$ all other $x_n$ can be found by using \eqref{1}.

\section{Absorption probabilities }
\label{s3}
\subsection{Absorption probabilities on a finite interval}
Consider our random walk on the finite interval $[a,b]$. Our main goal is the interval $[0,N]$, but we need more general results in the next sections.
Let $x_n=x_n^{[a,b]}$ be the expected number of visits to state $n \ \ (a\leq n\leq b).$
\begin{theorem} \label{t3} 
The probabilities of absorption in state $n$ ($a\leq n\leq b)$ in a $[pqrs]$ random walk on $[a,b]$, where $a$ and $b$ are absorbing barriers,  when starting in $i_0 \ \  (a< i_0< b)$ are\[s_nx_n\quad (s_a=s_b=1; s_i=s \ \  (i=a+1,\dots,b-1))\] where:
\begin{equation}\label{4}
x_a=x_a^{[a,b]}=\frac{\xi_1^{b-i_0}-\xi_2^{b-i_0}}{\xi_1^{b-a}-\xi_2^{b-a}}
\end{equation}
\begin{equation}\label{5}
x_b=x_b^{[a,b]}=\frac{\xi_2^{a-i_0}-\xi_1^{a-i_0}}{\xi_2^{a-b}-\xi_1^{a-b}}
\end{equation}
\begin{equation} \label{6}
x_n=x_n^{[a,b]}=\left\{\begin{array}{l}\frac{\zeta (\xi_2^{b-i_0}-\xi_1^{b-i_0})(\xi_1^a\xi_2^n-\xi_1^n\xi_2^a) }{\xi_1^b\xi_2^a-\xi_1^a\xi_2^b } \ \ \ \ \ \ \ \    (a+1\leq n \leq i_0) \\
\frac{\zeta (\xi_2^{a-i_0}-\xi_1^{a-i_0})(\xi_1^b\xi_2^n-\xi_1^n\xi_2^b)}{\xi_1^b\xi_2^a-\xi_1^a\xi_2^b} \ \ \ \ \ \ \ \    (i_0 \leq n\leq b-1)\end{array}\right.
\end{equation}
\end{theorem}
\begin{proof}
We start in $i_0 \ \ (a+2\leq i_0\leq b-2)$.
The set of difference equations:
\begin{equation*} 
(1 - r) x_ n  =\delta   (n,i_0 )+  px_{n- 1}+   qx_{ n+ 1} \quad	(a+2\leq n\leq b-2)			
\end{equation*}
has solutions (using Theorem \eqref{t1}):			
\[
x_n=\left\{\begin{array}{l}\zeta \xi_1^{n-i_0}+C_1 \xi_1^{n}+C_2 \xi_2^{n}\ \ \ \ \ \ \ \    (a+1\leq n \leq i_0) \\
\zeta \xi_2^{n-i_0}+C_1 \xi_1^{n}+C_2 \xi_2^{n}\ \ \ \ \ \ \ \    (i_0 \leq n\leq b-1)\end{array}\right.
\]		
Notice that the solution is also valid for $n=a+1$ and $n=b-1$, which can be seen by observing the difference equations for $n=a+2$ and $n=b-2: x_{a+1}$ and $x_{b-1}$ satisfy the difference pattern.
Using $(1-r)x_{a+1}=qx_{a+2}$ and $(1-r)x_{b-1}=px_{b-2}$ , where $b-2>a$, we get:
\[C_1=\frac{\zeta \xi_2^b(\xi_1^{a-i_0}-\xi_2^{a-i_0}) }{\xi_1^b\xi_2^a-\xi_1^a\xi_2^b}
\]
\[C_2=\frac{\zeta \xi_1^a(\xi_2^{b-i_0}-\xi_1^{b-i_0}) }{\xi_1^b\xi_2^a-\xi_1^a\xi_2^b}
\]
which leads to \eqref{6}.

$x_a=qx_{a+1}$ and $x_b=px_{b-1}$ leads to \eqref{4} and \eqref{5}.\\
After some calculations we find the result also valid for $i_o=a,a+1,b-1,b$.\\
For a state $j$ with absorption probability $s_j$ we have:\\
$P$(absorption in $j$ when starting in $i$)=$\sum_{k=0}^{\infty}p_{ij}^{(k)}s_j=s_jx_j.$\\
For an absorbing barrier we have $s_j=1$, so the probability of absorption in a barrier is $x_j$.
\end{proof}

After some calculations we have:
\[\sum_{n=a}^{b}s_nx_n=1\quad (s_a=s_b=1; s_i=s \ \  (i=a+1,\dots,b-1))\]

\subsection{Absorption probabilities on a semi-infinite interval}
We start with the semi-infinite interval $[a,\infty)$, where $a$ is an absorbing barrier. 
\begin{theorem}\label{t4}
The probability of absorption in a $[pqrs]$ random walk on $[a,\infty)$ when starting in $i_0$ is:\\
in absorbing barrier $a$:
\begin{equation}\label{extra}
\quad x_a^{[a,\infty)}=\xi_1^{a-i_0} \quad (a< i_0)
\end{equation}
in all other states $sx_n$, where
\begin{equation} \label{9}
x_n=x_n^{[a,\infty)}\left\{\begin{array}{l}\zeta \xi_1^{a-i_0}( \xi_1^{n-a}-\xi_2^{n-a})\ \ \ \ \ \ \ \    (a+1\leq n \leq i_0) \\
\zeta \xi_2^{n-a}( \xi_2^{a-i_0}-\xi_1^{a-i_0})\ \ \ \ \ \ \ \    (i_0 \leq n)\end{array}\right.
\end{equation}
\end{theorem}
\begin{proof}
We note: $P$(absorption in state $n$)=$sx_n\leq 1$ for fixed $s>0$. So $x_n$ is finite.
Consider a $[pqrs]$ random walk on $[a,\infty)$ where $a$ is an absorbing barrier.
We start in $i_0$ with $a+2\leq i_0.$
The set of difference equations:
\begin{equation*} 
(1 - r) x_ n  =\delta   (n,i_0 )+  px_{n- 1}+   qx_{ n+ 1} \quad	(a+2\leq n)			
\end{equation*}
has solutions (use Theorem \eqref{t1}, $\xi_1>1$ and $x_n$ is finite):			
\[
x_n=\left\{\begin{array}{l}\zeta \xi_1^{n-i_0}+C_2 \xi_2^{n}\ \ \ \ \ \ \ \    (a+1\leq n \leq i_0) \\
\zeta \xi_2^{n-i_0}+C_2 \xi_2^{n}\ \ \ \ \ \ \ \    (i_0 \leq n)\end{array}\right.
\]		
Notice that the solution is also valid for $n=a+1$ , which can be seen by observing the difference equations for $n=a+2$: $x_{a+1}$  satisfy the difference pattern.
Using $(1-r)x_{a+1}=qx_{a+2}$  we get: $C_2=-\zeta \xi_1^{a-i_0}\xi_2^{-a}$, which leads us to \eqref{9}.
$x_a=qx_{a+1}$ leads to \eqref{extra}.
\end{proof}
This result can also be obtained by taking $b\rightarrow \infty$ in \eqref{4} and \eqref{6}.\\ 
After some calculations we have:
\[\sum_{n=a}^{\infty}s_nx_n=1\quad (s_a=1; s_i=s \ (i>a))\]
In the next sections we also need absorbing probabilities on $(-\infty,b]$, where b is an absorbing barrier.
\begin{theorem} \label{t5}
The probability of absorption in the barrier $b$ in a $[pqrs]$ random walk on $(-\infty,b]$ when starting in $i_0$ is:
\begin{equation} \label{12}
 x_b^{(-\infty,b]}=\xi_2^{b-i_0} \quad (b> i_0)
 \end{equation}
 and probability of absorption in all other states is $sx_n$, where
 \[
x_n=x_n^{(-\infty,b]}=\left\{\begin{array}{l}\zeta \xi_1^{n-b}( \xi_1^{b-i_0}-\xi_2^{b-i_0})\ \ \ \ \ \ \ \    ( n \leq i_0) \\
\zeta \xi_2^{b-i_0}( \xi_2^{n-b}-\xi_1^{n-b})\ \ \ \ \ \ \ \    (i_0 \leq n \leq b-1)\end{array}\right.
\]
\end{theorem}
\begin{proof}
Proceed along the same lines as with ${[a,\infty)}$ 
\end{proof}

\subsection{Absorption probabilities on a infinite interval}
Our domain is the infinite interval $(-\infty,\infty).$
There are no absorption barriers.
\begin{theorem}
 Probability of absorption in state $n$ is $sx_n$, where
\begin{equation} \label{13}
x_n=x_n^{(-\infty,\infty)}=\left\{\begin{array}{l}\zeta \xi_1^{n-i_0}\ \ \ \ \ \ \ \    (n \leq i_0) \\
\zeta \xi_2^{n-i_0}\ \ \ \ \ \ \ \    (n \geq i_0)\end{array}\right.
\end{equation}	
\end{theorem}
\begin{proof}
Use Theorem \ref{t1} with $C_1=C_2=0$ (use: $\xi_1>1$, $0<\xi_2<1$ and $x_n$ is finite).
\end{proof} 
After some calculations we have:
\begin{equation}\label{dertien}\sum_{n=-\infty}^{\infty}sx_n=1
\end{equation}
\subsection{Asymptotic behavior of absorbing probabilities }
In this section we obtain asymptotic results for the probabilities of absorption  when $s \rightarrow 0$. We restrict to $(-\infty,\infty)$. The results for $[0,N]$ and $[0,\infty)$ are obtained by the same procedure.
\begin{lemma} \label{t11}
If $s \rightarrow 0$ then:\\
If $p>q:$\\  $\xi_1\sim \frac{p}{q}(1+\frac{s}{p-q}) \quad   \xi_2\sim 1-\frac{s}{p-q}+\frac{ps^2}{(p-q)^3} \quad \zeta \sim \frac{1}{p-q}[1-\frac{(p+q)s}{(p-q)^2}]$\\
If $p<q:\\ \xi_1\sim 1+\frac{s}{q-p} \quad   \xi_2\sim \frac {p}{q}[1-\frac{s}{q-p}+\frac{qs^2}{(q-p)^3}] \quad \zeta \sim \frac{1}{q-p}[1-\frac{(p+q)s}{(p-q)^2}]$\\ 
If $p=q$:\\
$\xi_1\sim 1+t+\frac{1}{2}t^2+\frac{1}{8}t^3$ \quad $ \xi_2\sim 1-t+\frac{1}{2}t^2-\frac{1}{8}t^3 \quad \zeta\sim \frac{1}{2pt}(1-\frac{1}{8}t^2)$ \\
where  $t=\sqrt{\frac{s}{p}}$.
\end{lemma}
\begin{proof}
We proof the last one:
$\zeta=[(1-r)^2-4 p^2]^{-\frac{1}{2}}=(4ps+s^2)^{-\frac{1}{2}}= (4ps)^{-\frac{1}{2}}(1+\frac{s}{4p})^{-\frac{1}{2}}
\sim \frac{1}{2\sqrt{ps}}(1-\frac{s}{8p})= \frac{1}{2pt}(1-\frac{1}{8}t^2)$\\
The rest goes in a similar way.
\end{proof}

\begin{theorem}
Asymptotic behavior of absorption probabilities on $(-\infty,\infty)$. We start in $i_0=0$.\\
 If $s \rightarrow 0$ and $p>q$:
\begin{equation*}
sx_{n}\sim  \frac{s}{p-q}\left\{1-\left [\frac{n}{p-q}+\frac{p+q}{(p-q)^2}\right ]s    \right\} \quad (n\geq 0)
\end{equation*}
\begin{equation*}
sx_{n}\sim  \frac{(\frac{p}{q})^ns}{p-q}\left\{1-\left [\frac{-n}{p-q}+\frac{p+q}{(p-q)^2}\right ]s    \right\} \quad (n\leq 0)
\end{equation*}
\begin{equation*}
\sum_{n=-1}^{-\infty}sx_n\sim \frac{qs}{(p-q)^2}\left\{1-\frac{2p+q}{(p-q)^2}s\right\}
\end{equation*}
\begin{equation}
\label{87}
\sum_{n=1}^{\infty}sx_n\sim 1-\frac{ps}{(p-q)^2}
\end{equation}
If $s \rightarrow 0$ and $p=q$:
\begin{equation*}
\sum_{n=1}^{\infty}sx_n=\sum_{n=-1}^{-\infty}sx_n\sim \frac{1}{2}-\frac{1}{4}\sqrt{\frac{s}{p}}
\end{equation*}
\begin{equation*}
sx_0\sim \frac{1}{2}\sqrt{\frac{s}{p}}
\end{equation*}
\end{theorem}
\begin{proof}
Without limitation we can take $i_0=0$.
We use  \eqref{13} and Lemma \ref{t11}. We proof \eqref{87}. In a similar way we can prove the other results. We have $p>q$, so:
$1-\xi_2\sim \frac{s}{p-q}[1-\frac{ps}{(p-q)^2}]$ , and
$(1-\xi_2)^{-1}\sim \frac{p-q}{s}[1+\frac{ps}{(p-q)^2}]$ and  $\sum_{n=1}^{\infty}sx_n=s\zeta \sum_{n=1}^{\infty}\xi_2^n=s\zeta\xi_2(1-\xi_2)^{-1}\sim$ \\
$ \frac{s}{p-q}[1-\frac{(p+q)}{(p-q)^2}s][1-\frac{s}{p-q}][\frac{p-q}{s}][1+\frac{ps}{(p-q)^2}]\sim 1-\frac{ps}{(p-q)^2}$
\end{proof}
\section{Expected time until absorption} \label{s4}
In this section we are interested in the expected time until absorption  in a $[pqrs]$ random walk with $pqs>0$. We define $m_i=m_i^D $ as the expected time until absorption  when starting in state $i$ on domain $D.$ In section    \ref{s3} we proved that absorption always occurs on domain $D$.
First time analysis gives:
\[m_i=p(m_{i+1}+1)+q(m_{i-1}+1)+r(m_i+1)+s.1\]
We start with a central Theorem.
\begin{theorem} \label{t7}
The set of difference equations
\begin{equation} \label{t6a}
(1-r)m_i=pm_{i+1}+qm_{i-1}+1\quad (i\in \mathbb Z)\quad (p+q+r+s=1, pqs>0)
\end{equation}
has solutions
\begin{equation}\label{t6b}
m_i=a\xi_1^{-i}+b\xi_2^{-i}+\frac{1}{s} \quad (i\in \mathbb Z)
\end{equation}
\end{theorem}
\begin{proof}
By substitution.
\end{proof}
The expected times until absorption are unique by the same argument given after Theorem \ref{t1}.
\begin{lemma}\label{lemma}
\[m_i^{[0,N]}\leq m_i^{[0,\infty)}\leq m_i^{(-\infty,\infty)}=\frac{1}{s}\]
\end{lemma}
\begin{proof}
Let $T_i^D$ be the time until absorption on domain $D$ when starting in $i \in D$. We have: $T_i^{[0,N]}\leq T_i^{[0,\infty]} \leq T_i^{(-\infty,\infty)}$, so $m_i^{[0,N]}\leq m_i^{[0,\infty)}\leq m_i^{(-\infty,\infty)}.$
We also have on $D=(-\infty,\infty)$:
 $m_i=m=s\sum_{k=1}^{\infty}k(1-s)^{k-1}=\frac{1}{s} $, which gives  the same result as:
$m=p(m+1)+q(m+1)+r(m+1)+s.1$
\end{proof}
\subsection{Expected time until absorption  on $[0,N]$}
\begin{theorem}\label{t8}
The expected time until absorption when starting in $i\ (i=0,1,\dots,N)$ in a $[pqrs]$ random walk on $[0,N]$ is:
\begin{equation}\label{15}
m_i=\frac{1}{s}\left \{\frac{(1-\xi_1^{N-i})(1-\xi_2^N)-(1-\xi_2^{N-i})(1-\xi_1^N)}{\xi_1^{N}-\xi_2^{N}}\right \}
\end{equation}
\end{theorem}
\begin{proof}
\[
(1-r)m_i=pm_{i+1}+qm_{i-1}+1\quad  (i=1,2,\dots,N-1)
\]
\[m_0=m_N=0\]
Use  theorem \ref{t7} to get the result.
\end{proof}

\subsection{Expected time until absorption  on $[0,\infty)$}
\begin{theorem} \label{t9}
The expected time until absorption when starting in $i\ (i=0,1,\dots)$ in a $[pqrs]$ random walk on $[0,\infty)$ is:
\begin{equation} \label{16}
m_i=\frac{1}{s}(1-\xi_1^{-i})
\end{equation}
\end{theorem}
\begin{proof}
\[
(1-r)m_i=pm_{i+1}+qm_{i-1}+1\quad  (i=1,2,\dots)
\]
\[m_0=0\]
Use Lemma \ref{lemma} and Theorem \ref{t7} (with $b=0$). 
\end{proof}
We get the same result by taking $N\rightarrow \infty$ in  theorem \ref{t8}.

\subsection{Asymptotic behavior of expected time until absorption}
In this section we obtain asymptotic results for the expected time until absorption  when $s \rightarrow 0$. We restrict to $[0,\infty)$. The results for $[0,N]$ and $(-\infty,\infty)$ are obtained by the same procedure.
$[0,\infty)$ is interesting because in this case we encounter both recurrent and transient results.

\begin{theorem}
Asymptotic behavior of expected time until absorption on $[0,\infty)$:\\
Case $p>q$:\\
 \[m_i\sim \frac{1-(\frac{q}{p})^i}{s}+\frac{i(\frac{q}{p})^i}{p-q}     \quad (s\rightarrow 0)\]
Case $p<q$:\\
 \[m_i\sim \frac{i}{q-p}  +\frac{1}{2}i\left\{\frac{p+q}{(q-p)^3  }-\frac{i}{(q-p)^2  } \right\} s     \quad (s\rightarrow 0)\]
Case $p=q$:
\[m_i \sim i\sqrt{\frac{1}{ps}}-\frac{i^2}{2p}    \quad (s\rightarrow 0)\]
\end{theorem}
\begin{proof}
Equation \eqref{16} states:
\begin{equation*} 
m_i=\frac{1}{s}(1-\xi_1^{-i})
\end{equation*}
We use Lemma \ref{t11}.\\
Case $p>q$:\\
$\xi_1^{-1}=\frac{q}{p}\xi_2\sim \frac{q}{p}[1-\frac{s}{p-q}]; \quad   \xi_1^{-i}\sim (\frac{q}{p})^i[1-\frac{is}{p-q}]$\\ 
Case $p<q$:\\
$\xi_1^{-1}=\frac{q}{p}\xi_2\sim  1-\frac{s}{q-p}+\frac{qs^2}{(q-p)^3} ; \quad    \xi_1^{-i}\sim 1- \frac{is}{q-p}  -\frac{1}{2}i\left\{\frac{p+q}{(q-p)^3  }-\frac{i}{(q-p)^2  } \right\}s^2      $ \\
where we used:
\begin{equation}
\label{76}
(1+as+bs^2)^i\sim 1+\binom{i}{1}(as)+\binom{i}{1}(bs^2)+\binom{i}{2}(as)^2 \quad (s\rightarrow 0)
\end{equation}
Case $p=q$:\\
$\xi_1^{-1}=\xi_2\sim 1- \sqrt{\frac{s}{p}}+\frac{s}{2p}; \quad  \xi_1^{-i}\sim 1-i\sqrt\frac{s}{p}+\frac{i^2s}{2p}  \    $
where we again used \eqref{76}.
\end{proof}
\section{Maximum and minimum of extended ruin problem} \label{s5}
Let $m$ be the minimum and $M$ be the maximum value of the random walk.

\subsection{Maximum and minimum on [0,N]}

\begin{theorem}
On the interval $[0,N]$ we have:
\[P(m=0)=\frac{\xi_1^{N-i_0}-\xi_2^{N-i_0}}{\xi_1^{N}-\xi_2^{N}}
\]
\[P(m=a)=(\xi_1^{N-i_0}-\xi_2^{N-i_0})\{\frac{1}{\xi_1^{N-a}-\xi_2^{N-a}}-\frac{1}{\xi_1^{N-a+1}-\xi_2^{N-a+1}}\} \quad (1\leq a\leq i_0)\]
\[P(M=b)=(\xi_2^{-i_0}-\xi_1^{-i_0})\{\frac{1}{\xi_2^{-b}-\xi_1^{-b}}-\frac{1}{\xi_2^{-b-1}-\xi_1^{-b-1}}\} \quad (i_0\leq b\leq N-1)\]
\[P(M=N)=\frac{\xi_2^{-i_0}-\xi_1^{-i_0}}{\xi_2^{-N}-\xi_1^{-N}}
\]
\end{theorem}
\begin{proof}
\[P(m=0)=x_0^{[0,N]}
\]
\[P(M=N)=x_N^{[0,N]}
\]
We notice: $\{m\leq a\}=\{$random walk visits $a$ after $n$ steps for some $n\geq 0 \}$, where $1\leq a\leq i_0$. 
If $a$ is not an absorbing barrier then we can detect a visit to $a$ by transforming $a$ in an absorbing barrier. 
We get:
\[P(m=a)=x_a^{[a,N]}-x_{a-1}^{[a-1,N]} \quad (1\leq a\leq i_0)\]
We can apply the same procedure in case of $M\geq b$. and get:
\[P(M=b)=x_b^{[0,b]}-x_{b+1}^{[0,b+1]} \quad (i_0\leq b\leq N-1)\]

\end{proof}
\subsection{Maximum and minimum on $[0,\infty)$} 
Using the same techniques as in the preceeding section, we get:
\begin{equation}\label{10}
P(m=0)=x_0^{[0,\infty)}=\xi_1^{-i_0}
\end{equation}
\[P(m=a)=x_a^{[a,\infty)}-x_{a-1}^{[a-1,\infty)}=\xi_1^{a-i_0-1}(\xi_1-1)
\quad (1\leq a\leq i_0) \]
\[P(M=b)=(\xi_2^{-i_0}-\xi_1^{-i_0})\{\frac{1}{\xi_2^{-b}-\xi_1^{-b}}-\frac{1}{\xi_2^{-b-1}-\xi_1^{-b-1}}\} \quad (b\geq i_0)\]
\subsection{Maximum and minimum on $(-\infty,\infty)$}
In a similar way we find:
$P(m=a)=x_a^{[a,\infty)}-x_{a-1}^{[a-1,\infty)}=\xi_1^{a-1-i_0}(\xi_1-1)\quad (a\leq i_0) $\\
$P(M=b)=x_b^{(-\infty,b]}-x_{b+1}^{(-\infty,b+1]}=\xi_2^{b-i_0}(1-\xi_2)\quad (b\geq i_0) $

\section{Value of the game}\label{s6}
Let $\mathfrak{n}$ be our final position after absorption. We define the value $v_i$ of the game as the expected value of $\mathfrak{n}$ when starting in $i$.
First step analysis shows:
\begin{equation}
\label{19}
(1-r)v_i=pv_{i+1}+qv_{i-1}+si
\end{equation}

\begin{theorem}
\begin{equation}
v_i=i+(p-q)m_i
\end{equation}
\end{theorem}
\begin{proof}
Substitute $v_i=i+(p-q)w_i$ in \eqref{19}. We obtain:\\
$(1-r)w_i=pw_{i+1}+qw_{i-1}+1$.\\
On interval $[a,b]$ we have: $v_a=a$ and $v_b=b$, so $w_a=w_b=0$. We get the same difference equations and conditions as for $m_i$, and so $w_i=m_i$ gives us a solution. The solution is unique (see remark after Theorem \ref{1}). \\
Taking $a=0$ and $ b  \rightarrow \infty$, we get a solution for $[0,\infty)$ (proof by substitution), which is also unique by the same arguments.\\
The same procedure works for $(-\infty,\infty)$.
\end{proof}

Let $g$ be the expectation  of the final gain.  $g_i=v_i-i=(p-q)m_i$ where $(p-q)$ is the unit time gain expectation  and $m_i$ is the expected duration of the game when starting in $i$.\\
The substitution $v_i=i+(p-q)w_i$ is suggested by the result of the calculation of $v_i$ (see Appendix A).

\section{A conjugate random walk} \label{s7}

We define $P_r=\frac{p}{1-r}, \ Q_r=\frac{q}{1-r},\ S_r=\frac{s}{1-r}.$ Besides  our original $[pqrs]$ random walk with $p+q+r+s=1$ and $pqs>0$ we also consider the conjugate $[P_rQ_rS_r]$ walk with $P_r+Q_r+S_r=1$ and $P_rQ_rS_r>0$.
We define $
\Xi_i(P_r,Q_r,S_r)=\xi_i(P_r,Q_r,0,S_r)\ (i=1,2)
$ and \\
$Z(P_r,Q_r,S_r)=\zeta(P_r,Q_r,0,S_r)$

 \begin{lemma} \label{t15}
$\Xi_i(P_r,Q_r,S_r)=\xi_i(p,q,r,s)\ (i=1,2) $ and $Z(P_r,Q_r,S_r)=(1-r)\zeta(p,q,r,s)$
\end{lemma}
\begin{proof}
Theorem \eqref{t1} gives:\\
$\xi_i(p,q,r,s)=\frac{(1-r)+(-1)^{i-1}[(1-r)^2-4pq]^{-\frac{1}{2}}}{2q}=\frac{1+(-1)^{i-1}[1-4P_rQ_r]^{-\frac{1}{2}}}{2Q_r}=$\\ $\xi_i(P_r,Q_r,0,S_r)\ (i=1,2).$ We also have:\\
$(1-r)\zeta(p,q,r,s)=(1-r)   [(1-r)^2-4pq]^{-\frac{1}{2}}=(1-4P_rQ_r)^{-\frac{1}{2}}=\zeta(P_r,Q_r,0,S_r)$
\end{proof}

\begin{theorem}
The $[pqrs]$ walk and the conjugate $[P_r,Q_r,S_r]$ walk gives the same results for: maximum and minimum of the walks, absorption probabilities, asymptotic behavior and the value of the game. The expected time until absorption in the conjugate case is $(1-r)$ times the expected time until absorption in the original walk.
\end{theorem}
\begin{proof}
By Lemma \ref{t15} we have: all results with only $\xi_i \ (i=1,2)$ in it will hold for both walks. For example: all formulas with relation to maximum and minimum. But there is more. Absorption probabilities are given by $sx_n$ and in \eqref{6} \eqref{9} \eqref{13} we see that these probabilities are always of the form $s\zeta F$ where $F$ is a function of $\xi_i \  (i=1,2)$. By Lemma \ref{t15} we have $s \zeta (p,q,r,s)=\frac{s}{1-r}Z(P_r,Q_r,S_r)= S_rZ(P_r,Q_r,S_r)  $, so the  $ s\zeta $ part in our original formulas can be changed to $S_rZ$ in the conjugate walk, which doesn't change the value. 
The section about asymptotic behavior also stays unchanged: $\frac{p}{q}=\frac{P_r}{Q_r}$ and $\frac{s}{p}=\frac{S_r}{P_r}$.
 The value of the conjugate game is the same as the value of the original game: $\frac{p-q}{s}=\frac{P_r-Q_r}{S_r}$.
The expected time until absorption needs some attention. The basis of all calculations in section \ref{s5} is Theorem \ref{t7}. Besides the $\xi_i \ (i=1,2)$ we have a term $\frac{1}{s}$ in the original walk. This will be changed in $\frac{1}{S}=\frac{1-r}{s}$ in the conjugate one, and all the formulas in the delayed walk are of the form $\frac{G}{s}$, where $G$ is a function of $\xi_i \ (i=1,2)$ so the expected time until absorption in the conjugate case is $(1-r)$ times the expected time until absorption in the original walk.
\end{proof}
\clearpage
\appendix
\section{Value of the game on $[0,N]$ }
\begin{theorem} \label{t12}
\begin{equation} \label{17}
v_{i_0}=i_0+\frac{(p-q)}{s}\left\{1-\frac{\xi_1^{N-i_0}(1-\xi_2^N)+\xi_2^{N-i_0}(\xi_1^N-1)}{\xi_1^N-\xi_2^N}\right \}
\end{equation}
\end{theorem}
\begin{proof}
Calculating $\sum_{n=1}^{i_0-1}\xi^n$ and $\sum_{n=i_0}^{N-1}\xi^n$
 and differentiating we get:
 
  \[\sum_{n=1}^{i_0-1}n\xi^{n-1}=\frac{1-\xi^{i_0}-i_0(1-\xi)\xi^{i_0-1}}{(1-\xi)^2}\]
   \[\sum_{n=i_0}^{N-1}n\xi^{n-1}=\frac{\xi^{i_0}-\xi^N-N(1-\xi)\xi^{N-1}+i_0(1-\xi)\xi^{i_0-1}}{(1-\xi)^2}\]
\begin{multline}
v_{i_0}=s\sum_{n=1}^{N-1}nx_n+Nx_N=s\sum_{n=1}^{i_0-1}\frac{n\zeta(\xi_2^{N-i_o}-\xi_1^{N-i_0})(\xi_2^{n}-\xi_1^{n})}{\xi_1^{N}-\xi_2^{N}}+ \\
s\sum_{n=i_0}^{N-1}\frac{n\zeta(\xi_2^{-i_o}-\xi_1^{-i_0})(\xi_1^{N}\xi_2^{n}-\xi_1^{n}\xi_2^{N})}{\xi_1^{N}-\xi_2^{N}}+\frac{N\xi_1^{N}\xi_2^{N}(\xi_2^{-i_o}-\xi_1^{-i_0})}{\xi_1^{N}-\xi_2^{N}}= \\
\frac{s \zeta (\xi_2^{N-i_o}-\xi_1^{N-i_0})}{\xi_1^{N}-\xi_2^{N}}
\left\{ \frac{\xi_2[1-\xi_2^{i_0}-i_0(1-\xi_2)\xi_2^{i_0-1}]}{(1-\xi_2)^2}+
\right.\\
- \left.
\frac{\xi_1[1-\xi_1^{i_0}-i_0(1-\xi_1)\xi_1^{i_0-1}]}{(1-\xi_1)^2}\right \} + \\
\frac{s \zeta(\xi_2^{-i_o}-\xi_1^{-i_0})}{\xi_1^{N}-\xi_2^{N}}
\left \{\xi_1^N \xi_2 \left [\frac{\xi_2^{i_0}-\xi_2^N-N(1-\xi_2)\xi_2^{N-1}+i_0(1-\xi_2)\xi_2^{i_0-1}}{(1-\xi_2)^2}\right ]+\right. \\
\left. -\xi_2^N \xi_1\left [ \frac{\xi_1^{i_0}-\xi_1^N-N(1-\xi_1)\xi_1^{N-1}+i_0(1-\xi_1)\xi_1^{i_0-1}}{(1-\xi_1)^2}\right ]\right \}
+ \\
\frac{N\xi_1^{N}\xi_2^{N}(\xi_2^{-i_o}-\xi_1^{-i_0})}{\xi_1^{N}-\xi_2^{N}}
\end{multline}

We first concentrate on the terms linear in $N$ : a  calculation shows that these vanishes. Next we concentrate on terms linear in $i_0$ : a calculation reduces to $i_0$. The remaining terms can be written as (after some calculation):
\[
\frac{s\zeta}{\xi_1^N-\xi_2^N}\left\{  \frac{\xi_2\Phi}{(1-\xi_2)^2}- \frac{\xi_1\Phi}{(1-\xi_1)^2}   \right\}=\frac{(p-q)\Phi}{s(\xi_1^N-\xi_2^N)}
\]
where $\Phi=(\xi_1^N-\xi_2^N)-(\xi_1^{N-i_0}-\xi_2^{N-i_0})+\xi_1^N\xi_2^N(\xi_1^{-i_0}-\xi_2^{-i_0})   $
\end{proof}

\end{document}